**Shifted Hankel determinants of Catalan numbers and related results**

Johann Cigler

## Abstract

In this (partly expository) paper we give a short overview about the close relationship between the sequence of Catalan numbers and Hankel determinants from the point of view of orthogonal polynomials and show that an analogous situation exists for more general sequences.

## 1. Introduction

The shifted Hankel determinants $\det\left(C_{n+i+j}\right)_{i,j=0}^{k-1}$ of the Catalan numbers $C_n = \frac{1}{n+1}\binom{2n}{n}$ satisfy (cf. [4], Th. 33)

$$\det\left(C_{n+i+j}\right)_{i,j=0}^{k-1} = \prod_{1\le i\le j\le n-1} \frac{2k+i+j}{i+j}. \tag{1.1}$$

For fixed $n$ the right-hand side of (1.1) is a polynomial $H_n(k)$.

Motivated by this result we look for other sequences whose shifted Hankel determinants are given by polynomials.

First we consider the numbers

$$M_b(n) = \sum_{k=0}^{n} \begin{bmatrix} n \\ k \end{bmatrix} b^{n-k}, \tag{1.2}$$

$b \in \mathbb{R}$, where the coefficients $\begin{bmatrix} n \\ k \end{bmatrix} = \binom{n+k}{k} - \binom{n+k}{k-1}$ are the entries of the Catalan triangle [7], A009766:

$$\left(\begin{bmatrix} n \\ k \end{bmatrix}\right)_{n,k\ge 0} = \begin{pmatrix} 1 & & & & & \\ 1 & 1 & & & & \\ 1 & 2 & 2 & & & \\ 1 & 3 & 5 & 5 & & \\ 1 & 4 & 9 & 14 & 14 & \\ \vdots & \vdots & \vdots & \vdots & \vdots & \ddots \end{pmatrix}. \tag{1.3}$$

For $b = 0$ this sequence reduces to the Catalan numbers $M_0(n) = C_n$.

For other values of $b$ we also get known sequences. For example $M_1(n)=C_{n+1}$, $M_2(n)=\binom{2n+1}{n}$. For $b=-1$ we get the Fine numbers [7], A000957, for $b=3$ the sequence [7], A049027 and for $b=4$ the sequence [7], A076025.

For these sequences there exist polynomials $H_n(b,x)$ such that

$$\det\left(M_b(n+i+j)\right)_{i,j=0}^{k-1}=H_n(b,k). \tag{1.4}$$

They are given by

$$H_n(b,x)=\det\left(\binom{x+i+j}{2j}+b\binom{x+i+j}{2j+1}\right)_{i,j=0}^{n-1}. \tag{1.5}$$

For $b=0$ we reproduce a combinatorial proof due to Christian Krattenthaler, which establishes a close connection with plane partitions of the form $(n-1,n-2,\cdots,1)$: There are bijections of such plane partitions with a) families of nonintersecting Dyck paths which are counted by the Hankel determinant (1.4) and b) families of other lattice paths which are counted by (1.5). The last property has previously been proved by Robert Proctor [8] using results of representation theory.

Finally we show similar results for central binomial coefficients in place of Catalan numbers. The analog of (1.1) is

$$\frac{\det\left(\binom{2n+2i+2j}{n+i+j}\right)_{i,j=0}^{k-1}}{2^{k-1}}=2^n\prod_{j=1}^{n-1}\frac{\binom{2k-1+2j}{j}}{\binom{2j}{j}}. \tag{1.6}$$

More generally for

$$M(b,n)=\sum_{j=0}^{n}\binom{n+j}{j}b^{n-j} \tag{1.7}$$

there are polynomials $V_n(b,x)$ such that

$$\frac{\det\left(M(b,n+i+j)\right)_{i,j=0}^{k-1}}{\det\left(M(b,i+j)\right)_{i,j=0}^{k-1}}=\frac{\det\left(M(b,n+i+j)\right)_{i,j=0}^{k-1}}{(2-b)^{k-1}}=V_n(b,k). \tag{1.8}$$

They are given by

$$V_n(b,x)=2^nH_n\left(b,\frac{2x-1}{2}\right)=\det\left(\binom{x+i+j}{2j}\frac{2x+2i}{x+i+j}+b\binom{x+i+j}{2j+1}\frac{2x+2i-1}{x+i+j}\right)_{i,j=0}^{n-1}.$$

I want to thank Sam Hopkins and Michael Somos for their helpful answers to my question [6] and Christian Krattenthaler for his combinatorial proof.

## 2. Some known background material.

**2.1.** A sequence of monic polynomials $p_n(x)$ of degree $n$ is called orthogonal if

$$p_n(x) = (x - s_{n-1})p_{n-1}(x) - t_{n-2}p_{n-2}(x) \tag{2.1}$$

for some values $s_n$ and $t_n$. Let $\Lambda$ denote the linear functional on the polynomials defined by $\Lambda\left(p_n(x)\right) = [n = 0]$. Then (2.1) implies the orthogonality relations $\Lambda\left(p_n(x)p_m(x)\right) = 0$ for $n \neq m$ and $\Lambda\left(p_n(x)^2\right) \neq 0$.

Define numbers $c(n,k)$ by

$$x^n = \sum_{k=0}^{n} c(n,k)p_k(x). \tag{2.2}$$

This implies

$$\begin{aligned} &c(0,k) = [k = 0], \\ &c(n,0) = s_0c(n-1,0) + t_0c(n-1,1), \\ &c(n,k) = c(n-1,k-1) + s_kc(n-1,k) + t_kc(n-1,k+1). \end{aligned} \tag{2.3}$$

The moments $M(n) = \Lambda\left(x^n\right)$ are given by

$$M(n) = \Lambda\left(x^n\right) = c(n,0). \tag{2.4}$$

The first two Hankel determinants are

$$\begin{aligned} &\det\left(M(i+j)\right)_{i,j=0}^{n-1} = \prod_{i=1}^{n-1}\prod_{j=0}^{i-1} t_j, \\ &\det\left(M(i+j+1)\right)_{i,j=0}^{n-1} = (-1)^n p_n(0)\prod_{i=1}^{n-1}\prod_{j=0}^{i-1} t_j. \end{aligned} \tag{2.5}$$

These results are well known. A detailed account can be found in [10] with a somewhat different notation.

**2.2.** The intimate relation between Hankel determinants and Catalan numbers stems from the case $s_n = 0$ and $t_n = 1$.

The corresponding monic orthogonal polynomials are the Fibonacci polynomials

$$F_n(x) = \sum_{j=0}^{\left\lfloor \frac{n}{2} \right\rfloor} (-1)^j \binom{n-j}{j} x^{n-2j}. \tag{2.6}$$

They satisfy

$$F_n(x) = xF_{n-1}(x) - F_{n-2}(x) \tag{2.7}$$

with initial values $F_0(x)=1$ and $F_1(x)=x.$

Since $c(n,k)=0$ if $k \not\equiv n \bmod 2$ we can write

$$x^n = \sum_{k=0}^{\left\lfloor \frac{n}{2} \right\rfloor} \left\langle {n \atop k} \right\rangle F_{n-2k} \tag{2.8}$$

with $\left\langle {n \atop k} \right\rangle = c(n,n-2k)$ if $0 \le 2k \le n$ and $\left\langle {n \atop k} \right\rangle = 0$ else.

These numbers satisfy the same recurrence

$$\left\langle {n \atop k} \right\rangle = \left\langle {n-1 \atop k-1} \right\rangle + \left\langle {n-1 \atop k} \right\rangle \tag{2.9}$$

as the binomial coefficients.

By induction it is easy to verify that

$$\left\langle {n \atop k} \right\rangle = c(n,n-2k) = \binom{n}{k} - \binom{n}{k-1} \tag{2.10}$$

for $2k \le n.$

As special case we get the moments

$$\Lambda\left(x^{2n}\right) = c(2n,0) = \left\langle {2n \atop n} \right\rangle = \binom{2n}{n} - \binom{2n}{n-1} = \frac{1}{n+1}\binom{2n}{n} = C_n. \tag{2.11}$$

The first terms of the table $\left(\left\langle {i \atop j} \right\rangle\right)$ (cf. [7], A008315) are

| | | | | | |
|---|---|---|---|---|---|
| | | | | | 1 |
| | | | | 1 | 0 |
| | | | 1 | 1 | |
| | | 1 | 2 | 0 | |
| | 1 | 3 | 2 | | |
| 1 | 4 | 5 | 0 | | |

Rows as printed (staggered triangle):

1

1 0

1 1

1 2 0

1 3 2

1 4 5 0

1 5 9 5

1 6 14 14 0

**2.3.** We also need the polynomials $f_n(x)=F_{2n}\left(\sqrt{x}\right)$ and $g_n(x)=\dfrac{F_{2n+1}\left(\sqrt{x}\right)}{\sqrt{x}}$.

The polynomials

$$f_n(x)=F_{2n}\left(\sqrt{x}\right)=\sum_{k=0}^{n}(-1)^k\binom{2n-k}{k}x^{n-k}=\sum_{j=0}^{n}(-1)^{n-j}\binom{n+j}{2j}x^j \tag{2.12}$$

satisfy

$$f_n(x)=(x-2)f_{n-1}(x)-f_{n-2}(x) \tag{2.13}$$

with initial values $f_0(x)=1$ and $f_1(x)=x-1$ and are orthogonal with $t_n=1$ and $s_0=1$ and $s_n=2$ for $n>0$.

From $x^{2n}=\sum_{k=0}^{n}c\left(2n,2n-2k\right)F_{2n-2k}(x)$ we get $x^n=\sum_{k=0}^{n}c\left(2n,2n-2k\right)f_{n-k}(x).$

If $\lambda$ denotes the linear functional defined by $\lambda\left(f_n\left(x\right)\right)=[n=0]$ then we get for the moments

$$M_0(n)=\lambda\left(x^n\right)=c(2n,0)=C_n. \tag{2.14}$$

The polynomials

$$g_n(x)=\frac{F_{2n+1}\left(\sqrt{x}\right)}{\sqrt{x}}=\sum_{k=0}^{n}(-1)^k\binom{2n+1-k}{k}x^{n-k}=\sum_{j=0}^{n}(-1)^{n-j}\binom{n+1+j}{2j+1}x^j \tag{2.15}$$

satisfy

$$g_n(x)=(x-2)g_{n-1}(x)-g_{n-2}(x) \tag{2.16}$$

with initial values $g_0(x)=1$ and $g_1(x)=x-2$ and are orthogonal with $t_n=1$ and $s_n=2$ for all $n\geq 0$.

From $x^{2n+1}=\sum_{k=0}^{n}c\left(2n+1,2n+1-2k\right)F_{2n+1-2k}(x)$ we get

$x^n=\sum_{k=0}^{n}c\left(2n+1,2n+1-2k\right)g_{n-k}(x).$ Therefore the moments of $g_n(x)$ are

$$M_1(n)=c(2n+1,1)=c(2n+2,0)=C_{n+1}. \tag{2.17}$$

From (2.5) we get the well-known result that $\det\left(C_{i+j}\right)_{i,j=0}^{n-1}=\det\left(C_{i+j+1}\right)_{i,j=0}^{n-1}=1.$ This characterizes the sequence of Catalan numbers as the uniquely determined sequence $a_n$ such that the Hankel determinants $\det\left(a_{i+j}\right)_{i,j=0}^{n-1}=\det\left(a_{i+j+1}\right)_{i,j=0}^{n-1}=1$ for all positive integers *n*.

## 3. The sequences $M_b(n)$.

**3.1.** Let us now consider more generally for $b \in \mathbb{R}$ the orthogonal polynomials $P_n(b,x)$ with $t_n = 1,\ \ s_0(b) = b+1$ and $s_n(b) = 2$ for $n > 0.$

Then we have $f_n(x) = P_n(0,x)$ and $g_n(x) = P_n(1,x).$

The polynomials $P_n(b,x)$ are monic and orthogonal and satisfy

$$P_n(b,x) = \left(x - s_{n-1}(b)\right) P_{n-1}(b,x) - P_{n-2}\left(b,x\right) \tag{2.18}$$

with $P_{-1}(b,x) = 0$ and $P_0(b,x) = 1.$

For $n > 0$ we get

$$P_n(b,x) = f_n(x) - b g_{n-1}(x). \tag{2.19}$$

For this holds for $n = 1$ because $P_1(b,x) = x - b - 1 = (x-1) - b = f_1(x) - b g_0(x).$ For $n > 1$ both sides satisfy the same recurrence.

From (2.12) and (2.15) we get

$$P_n(b,x) = f_n(x) - b g_{n-1}(x) = \sum_{j=0}^{n} (-1)^{n-j} \left( \binom{n+j}{2j} + b \binom{n+j}{2j+1} \right) x^j. \tag{2.20}$$

By (2.4) the moments of $P_n(b,x)$ are $M_b(n) = r(b,n,0),$ where

$$r(b,n,k) = r(b,n-1,k-1) + s_k r(b,n-1,k) + r(b,n-1,k+1) \tag{2.21}$$

with $r(b,n,k) = 0$ for $k < 0$ and $r(b,0,k) = [k = 0].$

From (2.21) it is easy to verify that

$$r(b,n,k) = \sum_{j=0}^{n-k} \binom{n+k+1+j}{j} \frac{n+k+1-j}{n+k+1+j} b^{n-k-j}. \tag{2.22}$$

This gives the moments

$$M_b(n) = \sum_{j=0}^{n} \binom{n+1+j}{j} \frac{n+1-j}{n+1+j} b^{n-j} = \sum_{j=0}^{n} \left( \binom{n+j}{j} - \binom{n+j}{j-1} \right) b^{n-j}. \tag{2.23}$$

The first terms are $1,\ 1+b,\ 2+2b+b^2,\ 5+5b+3b^2+b^3,\ 14+14b+9b^2+4b^3+b^4, \cdots.$

If we write $M_b(n) = \sum_{k=0}^{n} \begin{bmatrix} n \\ k \end{bmatrix} b^{n-k}$ then $\begin{bmatrix} n \\ k \end{bmatrix} = \binom{n+k}{k} - \binom{n+k}{k-1}.$

Note that $\begin{bmatrix} n \\ k \end{bmatrix} = \begin{bmatrix} n \\ k-1 \end{bmatrix} + \begin{bmatrix} n-1 \\ k \end{bmatrix}$ and $\sum_{k=0}^{n} \begin{bmatrix} n \\ k \end{bmatrix} = C_{n+1}.$

Since $t_n = 1$ we have $\det\left(M_b(i+j)\right)_{i,j=0}^{n-1} = 1.$

**3.2.** Let

$$C(x)=\frac{1-\sqrt{1-4x}}{2x}=\sum_{n\ge 0}C_n x^n \qquad (2.24)$$

be the generating function of the Catalan numbers which satisfies

$$C(x)=1+xC(x)^2. \qquad (2.25)$$

Then the generating series of the moments is

$$B_b(x)=\sum_{n\ge 0}M_b(n)x^n=\frac{C(x)}{1-bxC(x)}. \qquad (2.26)$$

Equation (2.21) is equivalent with

$$B_b(x)\left(C(x)-1)\right)^k=xB_b(x)\left(C(x)-1)\right)^{k-1}+2xB_b(x)\left(C(x)-1)\right)^k+xB_b(x)\left(C(x)-1)\right)^{k+1}$$

for $k>0$ and

$B_b(x)=1+bxB_b(x)+xB_b(x)\left(C(x)-1)\right)$ for $k=0$.

The first identity follows from

$$\left(C(x)-1)\right)=x+2x\left(C(x)-1)\right)+x\left(C(x)-1)\right)^2=x+2xC(x)-2x+xC(x)^2-2xC(x)+x$$
$$=xC(x)^2$$

which is true.

The second identity is also equivalent with

$$C(x)=\left(1-bxC(x)\right)+(b+1)xC(x)+xC(x)\left(C(x)-1\right)=1-bxC(x)+xC(x)+bxC(x)+xC(x)^2-xC(x)$$
$$=1+xC(x)^2.$$

From

$$\sum_{n\ge 0}M_2(n)x^n=\frac{C(x)}{1-2xC(x)}=\frac{1-\sqrt{1-4x}}{2x\sqrt{1-4x}}=\sum_{n\ge 0}\binom{2n+1}{n}x^n$$

we get $M_2(n)=\binom{2n+1}{n}$.

For $b=-1$ we get $\sum_{n\ge 0}M_{-1}(n)x^n=\frac{C(x)}{1+xC(x)}$. This implies that $M_{-1}(n)$ are the Fine numbers [7] A 000957.

## 4. Associated polynomials

**4.1.** The following well-known result (cf. e.g. [4], Th. 33)

$$\det\left(C_{n+i+j}\right)_{i,j=0}^{k-1}=\prod_{1\le i\le j\le n-1}\frac{2k+i+j}{i+j}=\prod_{j=1}^{n-1}\frac{\binom{2k+2j}{j}}{\binom{2j}{j}} \tag{3.1}$$

shows that the Hankel determinants on the left-hand side are polynomials in $k$ for fixed $n$.

The right-hand side can also be written as $\prod_{\ell=2}^{2n-2}\left(\frac{2k+\ell}{\ell}\right)^{\min\left(\left\lfloor\frac{\ell}{2}\right\rfloor,\left\lfloor\frac{2n-\ell}{2}\right\rfloor\right)}=\prod_{j=1}^{\left\lfloor\frac{n}{2}\right\rfloor}\prod_{i=2j}^{2n-2j}\frac{2k+i}{i}.$

Let us state this fact as

**Theorem 1**

*The polynomials*

$$H_n(x)=\prod_{j=1}^{n-1}\frac{\binom{2x+2j}{j}}{\binom{2j}{j}}=\prod_{1\le i\le j\le n-1}\frac{2x+i+j}{i+j}=\prod_{\ell=2}^{2n-2}\left(\frac{2x+\ell}{\ell}\right)^{\min\left(\left\lfloor\frac{\ell}{2}\right\rfloor,\left\lfloor\frac{2n-\ell}{2}\right\rfloor\right)}=\prod_{j=1}^{\left\lfloor\frac{n}{2}\right\rfloor}\prod_{i=2j}^{2n-2j}\frac{2x+i}{i} \tag{3.2}$$

*satisfy* $H_n(0)=1$ *and* $H_n(k)=\det\left(C_{n+i+j}\right)_{i,j=0}^{k-1}$ *for each positive integer* $k$.

The simplest proof of (3.1) uses Dodgson's condensation method (cf. [3], 2.3), which implies that $u(n,k)=\det\left(C_{n+i+j}\right)_{i,j=0}^{k-1}$ satisfies

$$u(n,k)u(n+2,k-2)-u(n+2,k-1)u(n,k-1)+u(n+1,k-1)^2=0. \tag{3.3}$$

We show first that $H_n(x)$ satisfies

$$H_n(x)H_{n+2}(x-2)-H_n(x-1)H_{n+2}(x-1)+H_{n+1}(x-1)^2=0. \tag{3.4}$$

We call (3.4) the *condensation property*.

$\frac{H_n(x)}{H_{n+1}(x-1)}=\frac{(2n)!}{n!}\frac{1}{\prod_{j=0}^{n-1}(2x+j)}$ and $\frac{H_n(x)}{H_{n+1}(x)}=\frac{(2n)!}{n!}\frac{1}{\prod_{j=n+1}^{2n}(2x+j)}$ give

$$\frac{H_n(x)}{H_{n+1}(x-1)}\frac{H_{n+2}(x-2)}{H_{n+1}(x-1)}=\frac{(2n)!}{n!}\frac{1}{\prod_{j=0}^{n-1}(2x+j)}\frac{(n+1)!}{(2n+2)!}\prod_{j=0}^{n}(2x+j-2)=\frac{(x-1)(2x-1)}{(2n+1)(2x+n-1)}$$

and

$$\frac{H_n(x-1)}{H_{n+1}(x-1)}\frac{H_{n+2}(x-1)}{H_{n+1}(x-1)}=\frac{(2n)!}{n!}\frac{1}{\prod_{j=n-1}^{2n-2}(2x+j)}\frac{(n+1)!}{(2n+2)!}\prod_{j=n}^{2n}(2x+j)=\frac{(x+n)(2x+2n-1)}{(2n+1)(2x+n-1)}$$

and therefore

$$\frac{H_n(x)}{H_{n+1}(x-1)}\frac{H_{n+2}(x-2)}{H_{n+1}(x-1)}-\frac{H_n(x-1)}{H_{n+1}(x-1)}\frac{H_{n+2}(x-1)}{H_{n+1}(x-1)}=1.$$

Now (3.1) follows with induction by $n$. For $n=0$ and $n=1$ (3.1) reduces to $\det\left(C_{i+j}\right)_{i,j=0}^{k-1}=1$ and $\det\left(C_{i+j+1}\right)_{i,j=0}^{k-1}=1$. Since $u(n+2,1)=C_{n+2}$ we can compute $u(n+2,k)$ for all positive integers $k$ using (3.3).

**4.2.** Another formula for $H_n(x)$ is given by

**Theorem 2**

$$H_n(x)=\det\left(\binom{x+i+j}{2j}\right)_{i,j=0}^{n-1}. \tag{3.5}$$

A direct proof has been given in [1], Th. 4. As observed there this is a special case of a more general result (cf. [5] and [9]).

**Lemma 3**

Let $p_n(x)=\sum_{j=0}^{n}(-1)^{n-j}p(n,j)x^j$ be monic orthogonal polynomials with moments $M(n)$. Then

$$\det\left(p(i+k,j)\right)_{i.j=0}^{n-1}=\frac{\det\left(M(n+i+j)\right)_{i,j=0}^{k-1}}{\det\left(M(i+j)\right)_{i,j=0}^{k-1}}. \tag{3.6}$$

(3.6) in combination with (2.12) gives

$$H_n(k):=\det\left(\binom{k+i+j}{2j}\right)_{i,j=0}^{n-1}. \tag{3.7}$$

**4.3. A combinatorial proof of (3.7)** (due to Christian Krattenthaler, personal communication).

**4.3.1.** Let us consider plane partitions of the form $(n-1, n-2, \cdots, 1)$ of integers between $0$ and $k$. These are arrays of the form

$$\begin{array}{llll} \pi_{1,1} & \pi_{1,2} & \cdots & \pi_{1,n-1} \\ \pi_{2,1} & \cdots & \pi_{2,n-2} & \\ \vdots & \iddots & & \\ \pi_{n-1,1} & & & \end{array}$$

with $0 \le \pi_{i,j} \le k$ such that $\pi_{i,j} \ge \pi_{i,j+1}$ and $\pi_{i,j} \ge \pi_{i+1,j}$.

For $k = 1$ the number of these partitions is $C_n$. For $n = 1$ we have only the empty partition.

For $n = 2$ we have $2 = C_2$ such partitions $\boxed{0}$, $\boxed{1}$.

For $n = 3$ we get $5 = C_3$ partitions: $\begin{array}{ll} \boxed{0} & \boxed{0} \\ \boxed{0} & \end{array}$, $\begin{array}{ll} \boxed{1} & \boxed{0} \\ \boxed{0} & \end{array}$, $\begin{array}{ll} \boxed{1} & \boxed{1} \\ \boxed{0} & \end{array}$, $\begin{array}{ll} \boxed{1} & \boxed{0} \\ \boxed{1} & \end{array}$, $\begin{array}{ll} \boxed{1} & \boxed{1} \\ \boxed{1} & \end{array}$.

As a typical example $(n = 6, k = 5)$ for the following proofs consider the partition

$$\begin{array}{lllll} 5 & 5 & 5 & 3 & 2 \\ 5 & 4 & 3 & 3 & \\ 2 & 2 & 2 & & \\ 2 & 1 & & & \\ 0 & & & & \end{array}$$

R. Proctor [8] using methods from representation theory has shown that the number of these plane partitions is $H_n(k)$.

We want to give a more elementary proof of this fact.

We first associate to each such plane partition a family of $k$ non-intersecting Dyck paths as suggested by Sam Hopkins in [6].

We sketch this map using our example.

We draw dividing lines between different entries of the partition as in the left figure. Then we pull them apart in direction $(-1,1)$ and add vertical steps at the beginning and horizontal steps at the end such that the new beginning and end are on the diagonal $y = x$ as shown in the right figure. Finally we turn the figure 45 degrees.

This gives $k$ non-overlapping Dyck paths with initial points $A_i = (-2i, 0)$ and end points $E_i = (2n + 2j, 0)$ for $0 \le i, j \le k-1$.

By the Lindström-Gessel-Viennot theorem [2] the number of families of non-intersecting paths from the set $\{A_i\}$ to the set $\{E_j\}$ is $\det\left(c\left(A_i,E_j\right)\right)_{i,j=0}^{k-1}$, where $c(A_i,E_j)$ denotes the number of Dyck paths from $A_i$ to $E_j$.

Since the number of Dyck paths from $(0,0)$ to $(2n,0)$ is given by the Catalan number $C_n$ we get

$$\det\left(c\left(A_i,E_j\right)\right)_{i,j=0}^{k-1}=\det\left(C_{n+i+j}\right)_{i,j=0}^{k-1}.$$

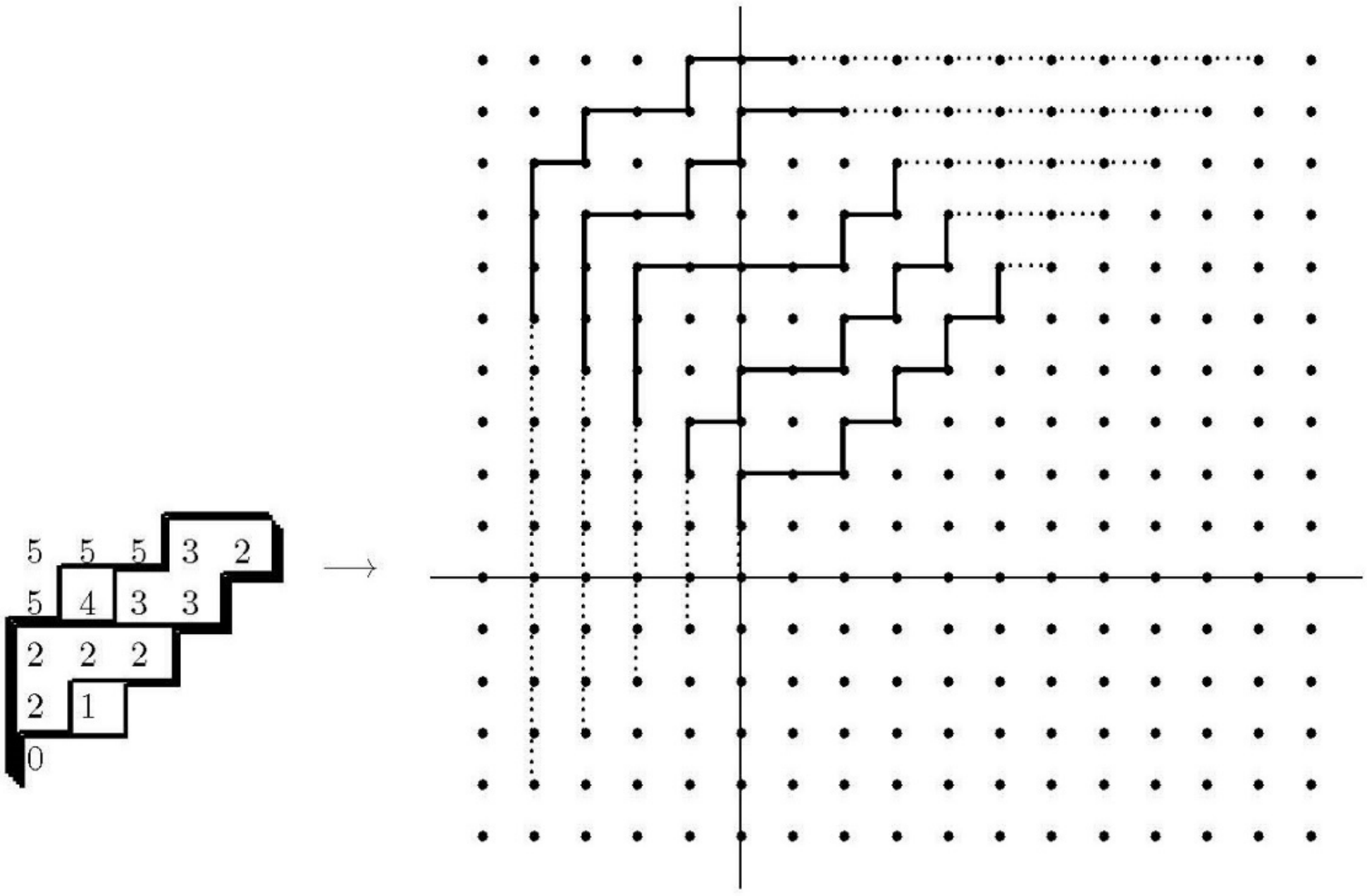

**4.3.2.** Next we show that the number of plane partitions is given by the determinant (3.7).

We assign to each row of the partition a lattice path which starts on the $y$-axis and ends on the $x$-axis, where the heights of the horizontal steps are the entries of the row.

If we denote by $H$ a horizontal step and by $V$ a vertical step then our example gives the paths (cf. the left side of the next figure)

$$H, HVHV, HHHVV, HVHVHHVVV, HHHVVHVHVV. \tag{3.8}$$

Then we shift these paths in direction $(2,1)$ and add horizontal and vertical steps at the start to obtain non-intersecting paths

$$P_i : A_i = (-1, k+i-1) \to E_i = (2i-1, i-1), \quad 1 \le i \le n,$$

and finally we can add a new vertical path $P_0$ from $(k-1,-1)$ to $(-1,-1)$ to obtain a nicer determinant (see the right part of the next figure).

In this way (3.8) is changed to

$$HVVVVVH,\ HHVVVHVHV, HHHVVVHHHVV,\\ HHHHHVHVHHVVV,\ HHHHHHHHVVHVHVV. \tag{3.9}$$

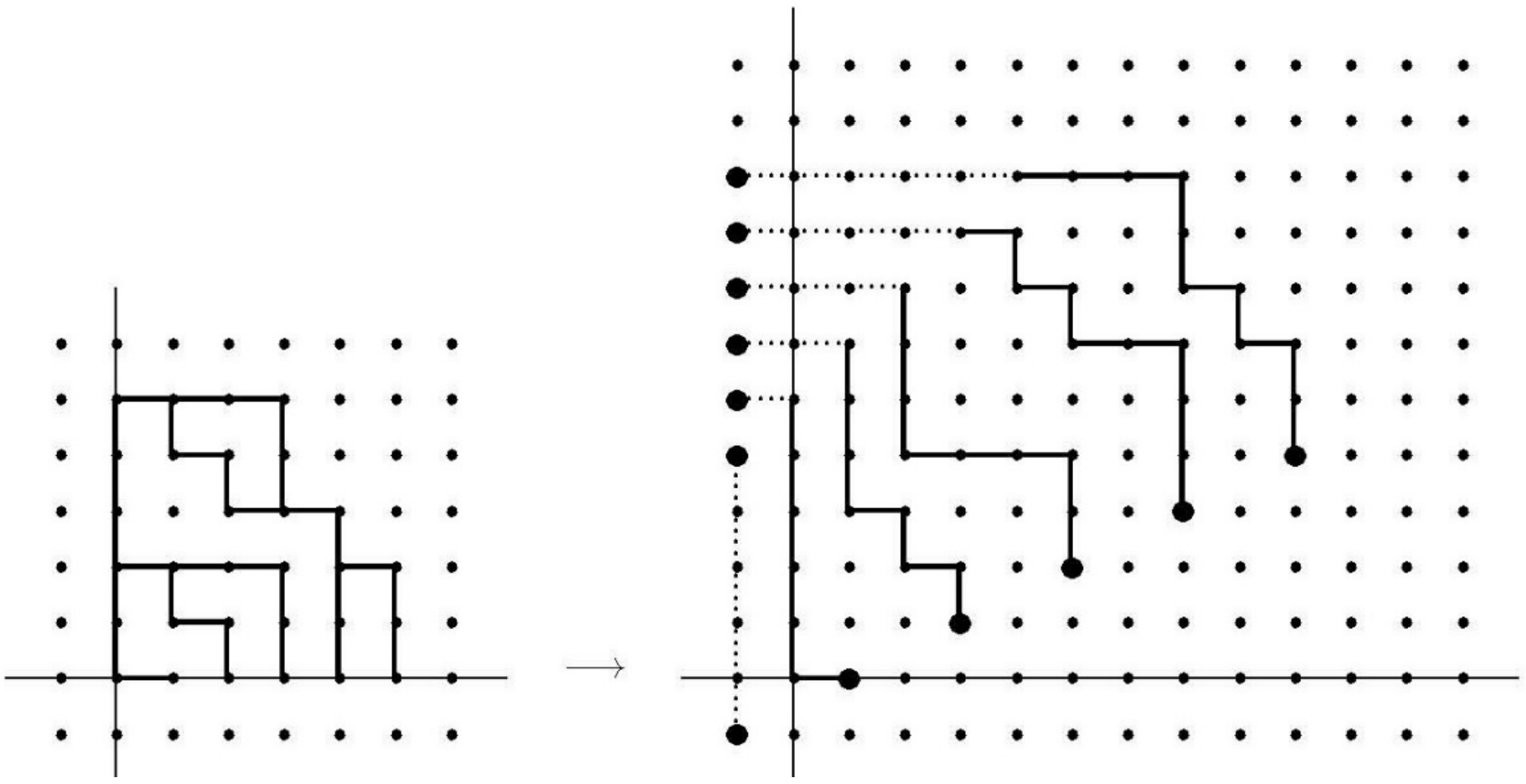

Let $s(A_i,E_j)$ denote the number of paths from $A_i$ to $E_j$ consisting of steps $H$ and $V$.

The Lindström-Gessel-Viennot theorem [2] says that

$\det\left(s\left(A_i,E_j\right)\right)_{i,j=0}^{n-1}$ is the number of all families of non-intersecting paths from the set $\{A_i\}$ to the set $\{E_j\}$. A path from $A_i=(-1,k+i-1)$ to $E_j=(2j-1,j-1)$ consists of $k+i+j$ steps from which $2j$ are horizontal. Therefore $s(A_i,E_j)=\binom{k+i+j}{2j}$.

This implies (3.7).

**4.4.** From (2.20) and (3.6) we get

**Theorem 4**

*The polynomials*

$$H_n(b,x)=\det\left(\binom{x+i+j}{2j}+b\binom{x+i+j}{2j+1}\right)_{i,j=0}^{n-1} \tag{3.10}$$

*satisfy* $H_n(b,k)=\det\left(H_{n+i+j}(b,1)\right)_{i,j=0}^{k-1}$ *with* $H_n(b,1)=M_b(n)$ *for a positive integer* $k$.

The first terms of the sequence $H_n(b,x)$ are $H_0(b,x)=1,\ H_1(b,x)=1+bx,$

$$H_2(b,x)=\frac{1}{6}(1+x)\left(6+b(b+6)x+2b^2x^2\right),$$

$$H_3(b,x)=\frac{(1+x)(2+x)\left(3+2x\right)}{180}\left(30+b(b^2+6b+30)x+3b^2(b+4)x^2+2b^3x^3\right).$$

For $b=0$ this reduces to Theorem 1 with $H_n(x)=H_n(0,x).$

For $b=1$ we have

$H_n(1,x)=H_{n+1}(0,x)$ because $H_n(1,k)=\det\left(C_{n+1+i+j}\right)_{i,j=0}^{k-1}=H_{n+1}(0,k)$ for all $k\in\mathbb{N}.$

**Theorem 5**

$$H_n(2,x)=2^nH_n\left(1,x-1/2\right)=\prod_{j=0}^{\left\lfloor\frac{n-1}{2}\right\rfloor}\prod_{i=2j+1}^{2n-2j-1}\frac{2x+i}{i}. \tag{3.11}$$

**Proof**

We know that both sides satisfy the condensation recursion. Let $h(n,x)=2^nH_n\left(1,x-1/2\right).$

We must show that $h(n,0)=\det\left(h(i+j,1)\right)_{i,j=0}^{k-1}=H_n(2,0),\ h(n,1)=\det\left(h(1+i+j,1)\right)_{i,j=0}^{k-1}$

For $n=0$ we have $H_0(2,x)=1$ and $H_n(2,1)=\binom{2n+1}{n}$ and $\det\left(H_{i+j}(2,1)\right)_{i,j=0}^{k-1}=1.$

For $n=1$ we have $H_1(2,x)=1+2x$ and $\det\left(H_{i+j+1}(2,1)\right)_{i,j=0}^{k-1}=1+2k.$

For $h(n,x)=2^nH_n\left(1,x-1/2\right)$ we have $h(0,x)=1$ and $h(1,x)=1+2x.$

$$h(n,1)=2^n\prod_{j=1}^{n}\frac{\binom{2j+1}{j}}{\binom{2j}{j}}=2^n\prod_{j=0}^{n}\frac{2j+1}{j+1}=\frac{(2n+1)!}{(n+1)!n!}=\binom{2n+1}{n}=H_n(2,1).$$

Changing $x\to\frac{2x-1}{2}$ gives $2^nH_n\left(1,x-1/2\right)=k_n\prod_{j=1}^{\left\lfloor\frac{n+1}{2}\right\rfloor}\prod_{i=2j}^{2n+2-2j}\left(2x+i-1\right)=k_n\prod_{j=0}^{\left\lfloor\frac{n-1}{2}\right\rfloor}\prod_{i=2j+1}^{2n-2j-1}(2x+i)$

for some $k_n$. Since $H_n(2,0)=1$ we get $H_n(2,x)=\prod_{j=0}^{\left\lfloor\frac{n-1}{2}\right\rfloor}\prod_{i=2j+1}^{2n-2j-1}\frac{2x+i}{i}.$

For $b\notin\{0,1,2\}$ the formulas are more complicated because then $H_n(b,x)$ has some irreducible factors of higher degree.

**4.5. Lemma 6 (cf. Michael Somos [6])**

*Let* $u(n,k)=\det\left(a_{n+i+j}\right)_{i,j=0}^{k-1}$. *Then condensation gives*

$$u(n,k)u(n+2,k-2)-u(n+2,k-1)u(n,k-1)+u(n+1,k-1)^2=0.$$

*If all* $u(n,k)\neq 0$ *then given* $u(0,k)$ *and* $u(1,k)$ *and* $u(n,1)$ *all* $u(n,k)$ *are uniquely determined.*

*Suppose the sequence of polynomials* $V_n(x)$ *satisfies the condensation condition*

$$V_n(x)V_{n+2}(x-2)-V_n(x-1)V_{n+2}(x-1)+V_{n+1}(x-1)^2=0$$

*and* $V_n(k)\neq 0$ *for all* $n,k.$

*Given* $a_0$ *and* $a_1$ *and* $V_0(x)$ *and* $V_1(x)$. *Then there are uniquely defined* $a_n$ *such that*

$$\det\left(a_{i+j}\right)_{i,j=0}^{k-1}=V_0(k) \text{ and } \det\left(a_{1+i+j}\right)_{i,j=0}^{k-1}=V_1(k).$$

**Corollary 7**

*Let* $V_n(b,x)=2^nH_n\left(b,\frac{2x-1}{2}\right)$. *Then*

$$V_n(b,k)=\frac{\det\left(V_{n+i+j}(b,1)\right)_{i,j=0}^{k-1}}{\det\left(V_{i+j}(b,1)\right)_{i,j=0}^{k-1}}. \tag{3.12}$$

Computations suggest that $V_n(b,1)=\sum_{j=0}^{n}\binom{n+j}{j}b^{n-j}$ and that these numbers are the moments corresponding to $s_b(0)=b+2,\ s_b(n)=2$ and $t_b(0)=2-b$ and $t_b(n)=1.$

We first prove

**Lemma 8**

*The orthogonal polynomials corresponding to* $s_b(0)=b+2,\ s_b(n)=2,\ t_b(0)=2-b$ *and* $t_b(n)=1$ *are*

$$p_n(b,x)=\sum_{j=0}^{n}(-1)^{n-j}\binom{n+j}{2j}\frac{2n}{n+j}x^j+b\sum_{j=0}^{n}(-1)^{n-j}\binom{n+j}{2j+1}\frac{2n-1}{n+j}x^j. \tag{3.13}$$

*Their moments are*

$$M(b,n)=\sum_{j=0}^{n}\binom{n+j}{j}b^{n-j} \tag{3.14}$$

*and their Hankel determinants are*

$$\det\left(M(b,i+j)\right)_{i,j=0}^{k-1}=(2-b)^{k-1}. \tag{3.15}$$

**Proof**

To prove (3.13) we must verify (2.1) which is easily done.

To prove (3.14) it suffices by (2.3) to compute $s(n,k)$ satisfying

$$s(0,k)=[k=0],$$
$$s(n,0)=(b+2)s(n-1,0)+(2-b)s(n-1,1),$$
$$s(n,k)=s(n-1,k-1)+2s(n-1,k)+s(n-1,k+1).$$

We have to show that $s(n,k)=\sum_{j=0}^{n-k}\binom{n+k+j}{j}b^{n-k-j}$ satisfies these conditions. This is easily verified.

$\det\left(M(b,i+j)\right)_{i,j=0}^{k-1}=(2-b)^{k-1}$ follows from (2.5).

If we write $\begin{bmatrix} n \\ j \end{bmatrix}=\binom{n+j}{j}$ then we get the triangle [7], A046899

$$\left(\begin{bmatrix} n \\ j \end{bmatrix}\right)_{n,j\ge 0}=\begin{pmatrix} 1 & & & & & \\ 1 & 2 & & & & \\ 1 & 3 & 6 & & & \\ 1 & 4 & 10 & 20 & & \\ 1 & 5 & 15 & 35 & 70 & \\ \vdots & \vdots & \vdots & \vdots & \vdots & \ddots \end{pmatrix}.$$

Note that $\begin{bmatrix} n \\ j \end{bmatrix}=\begin{bmatrix} n \\ j-1 \end{bmatrix}+\begin{bmatrix} n-1 \\ j \end{bmatrix}$ for $j<n$ and $\begin{bmatrix} n \\ n \end{bmatrix}=2\begin{bmatrix} n \\ n-1 \end{bmatrix}$.

This implies that $\sum_{j=0}^{k}\begin{bmatrix} n \\ j \end{bmatrix}=\begin{bmatrix} n+1 \\ k \end{bmatrix}$ and

$$M(1,n)=\sum_{j=0}^{n}\begin{bmatrix} n \\ j \end{bmatrix}=\begin{bmatrix} n+1 \\ n-1 \end{bmatrix}+\begin{bmatrix} n \\ n \end{bmatrix}=\binom{2n}{n-1}+\binom{2n}{n}=\binom{2n+1}{n}.$$

If we set $W_n(x)=\frac{V_n(x)}{2-b}$ then $V_n(k)=\frac{\det\left(V_{n+i+j}(1)\right)_{i,j=0}^{k-1}}{\det\left(V_{i+j}(1)\right)_{i,j=0}^{k-1}}=\frac{\det\left(V_{n+i+j}(1)\right)_{i,j=0}^{k-1}}{(2-b)^{k-1}}$

reduces to $W_n(k)=\frac{\det\left(V_{n+i+j}(1)\right)_{i,j=0}^{k-1}}{(2-b)^{k}}=\det\left(\frac{V_{n+i+j}(1)}{2-b}\right)_{i,j=0}^{k-1}=\det\left(W_{n+i+j}(1)\right)_{i,j=0}^{k-1}.$

**Proof of Corollary 7**

The conditions of Lemma 6 are satisfied:

1) $W_n(x)$ satisfies the condensation condition.

2) $W_0(x)=\frac{1}{2-b}$ and $W_1(x)=\frac{2+(2x-1)b}{2-b}.$

3) $\det\left(W_{i+j}(1)\right)_{i,j=0}^{k-1}=\frac{1}{2-b}$ and $\det\left(W_{i+j+1}(1)\right)_{i,j=0}^{k-1}=\frac{2+(2k-1)b}{2-b}$ by (2.5) and $(-1)^k p_k(b,0)=2+(2k-1)b.$

**Remark**

For $b=2$ we get from (3.11)

$$H_n\left(2,\frac{2x-1}{2}\right)=2^n H_n(1,x-1) \tag{3.16}$$

and therefore $2^n H_n\left(2,\frac{1}{2}\right)=\sum_{j=0}^{n}\binom{n+j}{j}2^{n-j}=4^n.$

From (3.6) we get

**Theorem 9**

$$V_n(b,x)=2^n H_n\left(b,\frac{2x-1}{2}\right)=\det\left(\binom{x+i+j}{2j}\frac{2x+2i}{x+i+j}+b\binom{x+i+j}{2j+1}\frac{2x+2i-1}{x+i+j}\right)_{i,j=0}^{n-1}$$

(3.17)

**Remark**

The last results are another illustration of the fact that there are many similarities between Fibonacci polynomials and Catalan numbers on the one hand and Lucas polynomials and central binomial coefficients on the other hand.

The analogs of the Fibonacci polynomials (2.6) are the (modified) Lucas polynomials defined by

$$L_n(x)=xL_{n-1}(x)-t_nL_{n-2}(x) \tag{3.18}$$

with initial values $L_0(x)=1$ and $L_1(x)=x$ where $t_0=2$ and $t_n=1$ for $n>0.$ The analog of (2.8) is

$$\sum_{k=0}^{\left\lfloor \frac{n}{2} \right\rfloor} \binom{n}{k} L_{n-2k}(x) = x^{n}. \tag{3.19}$$

If we define the linear functional $\Lambda$ on the polynomials by $\Lambda\left(L_n(x)\right) = [n=0]$ this implies

$$\Lambda\left(x^{2n}\right) = \binom{2n}{n} \tag{3.20}$$

and $\Lambda\left(x^{2n+1}\right) = 0.$

From

$$L_n(x) = \sum_{k=0}^{\left\lfloor \frac{n}{2} \right\rfloor} \binom{n-k}{k} \frac{n}{n-k} (-1)^k x^{n-2k} \tag{3.21}$$

for $n > 0$ we get

$$L_{2n}\left(\sqrt{x}\right) = \sum_{j=0}^{n} (-1)^{n-j} \binom{n+j}{2j} \frac{2n}{n+j} x^j \tag{3.22}$$

and

$$\frac{L_{2n-1}\left(\sqrt{x}\right)}{\sqrt{x}} = \sum_{j=0}^{n} (-1)^{n-1-j} \binom{n+j}{2j+1} \frac{2n-1}{n+j} x^j. \tag{3.23}$$

Thus (3.13) can be written as

$$p_n(b,x) = L_{2n}\left(\sqrt{x}\right) - b \frac{L_{2n-1}\left(\sqrt{x}\right)}{\sqrt{x}}, \tag{3.24}$$

which is an analog of (2.20).

## 5. A related result

**Theorem 10**

*Let*

$$h_n(x) = \prod_{k=1}^{n-1} \left( \frac{x+k}{k} \right)^{\min(k,n-k)}. \tag{4.1}$$

*Then we get* $h_n(1) = b(n) = \binom{n}{\left\lfloor \frac{n}{2} \right\rfloor}$ *and* $h_n(k) = (-1)^{n\binom{k}{2}} \det\left(b(n+i+j)\right)_{i,j=0}^{k-1}$. (4.2)

**Proof**

The polynomials $h_n(x)$ satisfy

$$(-1)^n h_n(x)h_{n+2}(x-2)-h_{n+2}(x-1)h_n(x-1)+h_{n+1}(x-1)^2=0.$$

Therefore $u(n,k)=(-1)^{n\binom{k}{2}}h_n(k)$ satisfies

$$u(n,k)u(n+2,k-2)-u(n+1,k-1)u(n,k-1)+u(n+1,k-1)^2=0.$$

It would be interesting if there are analogous results as above.

## References

[1] Johann Cigler, Some observations about determinants which are related with Catalan numbers and related topics, arXiv:1902.10468

[2] Ira M. Gessel and Gérard Viennot, Determinants, Paths and Plane Partitions, (1989 preprint), http://people.brandeis.edu/~gessel/homepage/papers/pp.pdf

[3] Christian Krattenthaler, Advanced determinant calculus, Sém. Loth. Comb. 42 (1999)

[4] Christian Krattenthaler, Advanced determinant calculus: A complement, Linear Algebra Appl. 411 (2005), 68-166

[5] Christian Krattenthaler, Hankel determinants of linear combinations of moments of orthogonal polynomials II, preprint

[6] MO, https://mathoverflow.net/questions/389657/some-nice-polynomials-related-to-hankel-determinants

[7] OEIS, The Online Encyclopedia of Integer Sequences, http://oeis.org/

[8] Robert A. Proctor, New Symmetric Plane Partition Identities from Invariant Theory Work of De Concini and Procesi, Europ. J. Comb. 11(1990), 289-300

[9] Mike Tyson, A theorem on shifted Hankel determinants, arXiv:1905.00327

[10] Gérard Viennot, Une théorie combinatoire des polynômes orthogonaux généraux, http://www.xavierviennot.org/xavier/polynomes_orthogonaux.html

Johann Cigler, Fakultät für Mathematik, Universität Wien, johann.cigler@univie.ac.at.